\def\draft{n}

\documentclass[12pt]{amsart}
\usepackage{fullpage,amssymb,epic,eepic,amscd,epsfig}

\theoremstyle{plain}

\newtheorem{theorem}{Theorem}
\newtheorem{proposition}{Proposition}[section]
\newtheorem{lemma}[proposition]{Lemma}

\theoremstyle{definition}

\theoremstyle{remark}

\newtheorem{remark}[proposition]{Remark}

\def\printname#1{
	\if\draft y
		\smash{\makebox[0pt]{\hspace{-0.5in}
			\raisebox{8pt}{\tt\tiny #1}}}
	\fi
}

\newcommand{\psdraw}[2]
        {\begin{array}{c} \hspace{-1.3mm}
	\raisebox{-4pt}{\psfig{figure=draws/#1.ps,width=#2}}
	\hspace{-1.9mm}\end{array}}
\newlength{\standardunitlength}
\catcode`\@=11
\long\def\@makecaption#1#2{%
    \vskip 10pt
    
\setbox\@tempboxa\hbox{
      \small\sf{\bfcaptionfont #1. }\ignorespaces #2}%
    \ifdim \wd\@tempboxa >\captionwidth {%
        \rightskip=\@captionmargin\leftskip=\@captionmargin
        \unhbox\@tempboxa\par}%
      \else
        \hbox to\hsize{\hfil\box\@tempboxa\hfil}%
    \fi}
\font\bfcaptionfont=cmssbx10 scaled \magstephalf
\newdimen\@captionmargin\@captionmargin=2\parindent
\newdimen\captionwidth\captionwidth=\hsize
\catcode`\@=12

\def\lbl#1{\label{#1}\printname{#1}}

\def\eqdef{\overset{\text{def}}{=}}

\def\eqdef{\overset{\text{def}}{=}}

\def\BZ{\mathbb Z}
\def\BQ{\mathbb Q}

\def\B{\mathcal B}

\def\K{\mathcal K}

\def\G{\mathcal G}
\def\O{\mathcal U}

\def\M{\mathcal M}

\def\ihs{integral homology 3-sphere}          

\def\fti{finite type invariant}

\def\la{\langle}
\def\ra{\rangle}

\def\FO#1{\mathcal F_{#1}\mathcal O}
\def\FV#1{\mathcal F_{#1}\mathcal V}

\def\GO#1{\mathcal G_{#1}\mathcal O}
\def\GV#1{\mathcal G_{#1}\mathcal V}

\def\cc{Chinese character}

\def\Ao{\mathcal A(\emptyset)}
\def\Ai{\mathcal A(S^1)}
\def\Bi{\mathcal B'}
\def\Bii{\mathcal B''}
\def\Bwh{\mathcal B_{wh}}
\def\reftool{\eqref{eq.wac}}

\def\v8{\vskip 8pt}
\def\lbl#1{\label{#1}}
\def\ZL{{ Z}^{ L  M  O}}
\def\ZK{{ Z}^{ K}}
\def\clos#1{\la #1 \ra}

\def\Ae#1{{\mathcal A}_{#1}(\emptyset)}

\def\Aps#1{{\mathcal A}'_{#1}(S^1)}
\def\Ap{{\mathcal A}'(S^1)}

\def\hpsi{\Psi}
\def\i{\iota}

\def\U{\mathcal U}

\def\ZL{ Z^{LMO}}
\def\ZK{Z^K}

\def\hb{\hbar}

\begin{document}


\title[The Alexander polynomial and finite type 3-manifold invariants]{
The Alexander polynomial and finite type 3-manifold invariants}

\author{Stavros Garoufalidis}
\address{Department of Mathematics \\
         Brandeis University \\
         Waltham, MA 02254-9110, U.S.A. }
\email{stavros@oscar.math.brandeis.edu}
\thanks{The  authors were partially supported by NSF grant 
       DMS-95-05105 and by the CNRS respectively. \newline
This and related preprints can also be obtained at
{\tt http://www.math.brown.edu/$\sim$stavrosg } and at
{\tt http://www.math.sciences.univ-nantes.fr/preprints/ }}

\author{Nathan Habegger}
\address{UMR 6629 du CNRS, Universit\'e de Nantes \\
         D\'epartement de Math\'e\-matiques \\
         2 rue de la Houssini\`ere \\
         44072 NANTES Cedex 03, France }
\email{habegger@math.univ-nantes.fr}

\date{Prliminary version, July 22,1997.  This version \today}

\begin{abstract}
Using elementary counting methods, we calculate the universal perturbative 
invariant (also known as the $LMO$ invariant) of a
3-manifold $M$, satisfying $H_1(M,\BZ)=\BZ$, in terms of the Alexander 
polynomial
of $M$. We show that $+1$ surgery on a knot in the 3-sphere induces an 
injective
map from finite type invariants of \ihs s to finite type invariants of knots. 
We also show that weight systems of degree $2m$ on knots, obtained by 
applying finite
type $3m$ invariants  of \ihs s, lie in the algebra of Alexander-Conway weight
systems, thus answering the questions raised in \cite{Ga}. 

\end{abstract}

\maketitle

\tableofcontents


\section{Introduction}
\lbl{sec.intro}

\subsection{History}
\lbl{sub.his}

In their fundamental paper, T.T.Q. Le, J. Murakami and T. Ohtsuki
\cite{LMO} constructed a map $\ZL$ which associates to every oriented 
3-manifold an element of the graded (completed) Hopf algebra $\Ao$ of trivalent 
graphs.\footnote{For a different construction of $\ZL(M)$ for a rational
homology 3-sphere $M$, see \cite{BGRT2}.}
The restriction of this map to the set of oriented \ihs s was shown in
\cite{L} to be the {\em universal} \fti\ of \ihs s 
(i.e., it classifies such invariants).  Thus $\ZL$ is a rich 
(though not fully understood) invariant of \ihs s. However, the invariant
$\ZL$ behaves differently as soon as the
first Betti number of the 3-manifold, $b_1(M)$, is positive. In \cite{Ha2},
the second author used an elementary counting argument to deduce that
$\ZL(M)=1$, if $b_1(M) > 3$, and to compute $\ZL(M)$, if
$b_1(M)=3$ (and also for $b_1(M)=2$, see \cite{BH}), in terms of the 
Lescop
invariant \cite{Le} of $M$.  
It is an open problem to compute $\ZL(M)$, for  $b_1(M) =0,1$.

It is the purpose of the present paper to exploit elementary
counting methods in order to calculate $\ZL(M)$, for 3-manifolds $M$ which
satisfy $H_1(M,\BZ)=\BZ$, in terms of a ``classical invariant'' of $M$, namely 
its
Alexander polynomial.  This includes the special case of 
0-{\em surgery}\footnote{given a framed link $L$ in a 3-manifold
$M$, we denote by $M_L$ the result of Dehn surgery on $L$.} of a
knot $K$ in $S^3$, $S^3_{K,0}$,
in which case the  Alexander polynomial of
$S^3_{K,0}$ is the Alexander-Conway polynomial of 
$K$.\footnote {An earlier version of 
this paper contained only this special
case.  We extend special thanks to C. Lescop, for help in extending to the
general case and to D. Thurston, for pointing out that the result should hold in 
this
generality.} An important ingredient of our computation is the recent result of 
A.
Kricker, B. Spence, and I. Aitchinson, \cite{Kr,KSA}, calculating the Conway 
weight
system on Chinese characters.

Although the invariant $\ZL(S^3_{K_,+1})$, of +1-surgery on a knot $K$ 
(in contrast to 0-surgery), is not
determined by the Alexander-Conway polynomial of $K$ (there are examples with
nontrivial invariant, and trivial Alexander polynomial), we show that after
truncating
$\ZL$ at degree $m$, the associated degree $2m$ knot {\em weight system}
lies in the algebra of the Alexander-Conway weight systems.
Similar methods allow us to show that \fti s of \ihs s are determined by their
associated knot invariants, thus answering positively the questions (see below)
that were posed in \cite{Ga} prior to the 
construction of
the LMO invariant.  (At that time, the only known finite type invariant of 
3-manifolds was the Casson invariant.)

\subsection{Statement of the results}
\lbl{sub.res}

All 3-manifolds and links considered in the present paper will be oriented.

\begin{theorem}
\lbl{thm.al}
Let $M$ be an oriented, closed, connected 3-manifold 
satisfying $H_1(M,\BZ)=\BZ$.
The universal invariant $\ZL(M) \in \Ao$ can be calculated in
terms of the Alexander-Conway polynomial $A(M)$ of the 3-manifold.
Conversely, the Alexander polynomial of $M$ can be
calculated in terms of the universal invariant $\ZL(M) \in \Ao$.
\end{theorem}
A precise formula relating the two invariants will be given in section  
\ref{sec.pre}.
\v8

We outline here the basic idea of the proof, which though somewhat 
technical, really is quite elementary:  If a manifold $M$ is obtained 
by $0$-surgery on a knot $K$ in $S^3$ (the general case of a manifold 
satisfying $H_1(M,\BZ)=\BZ$ is not much harder), then quite 
immediately from
the definitions, the degree $m$ part of $\ZL(M)$ can be computed from the part of the 
Kontsevich
integral of $K$ (written in Chinese characters) which has $2m$ legs and
$2m$ internal vertices.  Since there are no components which are 
intervals  (because of the
0-framing), and since by the anti-symmetry relation, all trees vanish, the 
only
contributing part consists of wheels.  But this part is known to determine the 
Alexander polynomial of $K$, and thus that of $M$.
\v8

Before we state the next result, we need to recall some standard definitions
and notation from the theory of \fti s of knots and \ihs s, see \cite{B-N,Oh,
Ga, LMO}.

Let $\K$ denote the vector space over $\BQ$ on the set of isotopy classes of
oriented knots in $S^3$ and let $\FV m$ (for a nonnegative integer $m$) denote
the vector space of finite type (i.e., Vassiliev)  invariants of knots of type
$m$, \cite{B-N}.   
Similarly, let $\M$ denote the vector space over $\BQ$ on the set
of orientation preserving diffeomorphism  classes of oriented \ihs s, and let
$\FO m$  denote the vector space of finite type (i.e., Ohtsuki) invariants of
\ihs s of type $m$, \cite{Oh}.   
In \cite{Ga} we considered the map $K \mapsto S^3_{K,+1}$. 
This is a classical map, often used in the study of knots (or 3-manifolds).
This yields a map $ \K \to \M$ and a dual map $\Phi: \M^\ast \to
\K^\ast$ (where $V^\ast$ denotes the dual of a vector space $V$).  In
\cite{Ga} the following questions were posed:

\v8
\begin{description}
\item[Q1]
Does the above map send
$\FO {3m}$ to $\FV {2m}$? 
\item[Q2]
Is the restriction of the map
$\Phi$ to $\FO {3m}$ one-to-one, for all $m$? 
\item[Q3]
Assuming the answer to Question 1 is affirmative, 
and given $v \in \FO {3m}$, is 
it true that
the associated degree
$2m$ knot weight system
lies in the algebra of the Alexander-Conway weight systems?
\end{description}

Let $v$ be a $\BQ$-valued
invariant of \ihs s and let $\Phi(v)$ be the associated invariant of knots
in $S^3$.  Question 1 asks whether $\Phi(v)$ is a finite type invariant of
knots in $S^3$ (together with an estimate of the type of the invariant), if $v$
is a finite type invariant of \ihs s. Question 2 asks whether $\Phi(v)$
determines $v$. (It should be noted, however, that there are \ihs s that
cannot be obtained by
$\pm 1$ surgery on a knot in $S^3$, see \cite{A}.)
Question 3 is concerned with the finite type knot invariant $\Phi(v)$ and
asks whether in degree $2m$ (the maximum possible degree by question 1), $\Phi(v)$ is a classical
knot invariant (on elements in the $2m$-th term of the Vassiliev filtration),
given by a polynomial in the Alexander-Conway coefficients.

Building on work of the first author and J. Levine (a preliminary version of 
\cite{GL2}),
Question 1 was answered affirmatively by the second author, \cite{Ha}.  
Alternative
proofs were later given in \cite{GL2,L}.  The methods used in
\cite{Ha} and
\cite{GL2} were a mixture of geometric topology together with a counting
argument. On the other hand, \cite{L} used the $\ZL$ invariant
and an elementary counting argument.  

Using elementary counting arguments similar to those in \cite{L}, together
with properties of the $\ZL$ invariant, enables us to show that
Questions 2 and 3 above are true.

\begin{theorem}
\lbl{thm.q2}
The association, which takes a knot in $S^3$ to the integral homology sphere
obtained by +1-framed surgery on the knot, induces an
injection from the space of finite type 3-manifold invariants (in the
sense of Ohtsuki), to the space of finite type (Vassiliev) knot
invariants. 
\end{theorem}

\begin{theorem}
\lbl{thm.q3}
Let $v$ be a finite type $3m$ invariant (in the sense of Ohtsuki) of homology
3-spheres.  Then the associated degree $2m$ knot weight system lies in the
algebra of the Alexander-Conway weight systems.
\end{theorem} 
 
\begin{remark}
\lbl{rem.gr}
Theorem  \ref{thm.q2} does not hold at the graded level, i.e.,
the associated graded map
$\G_m\Phi: \GO {3m}\to \GV {2m}$ is not one-to-one for $ m \geq 4$ (see remark 
\ref{rem.gr2}).
\end{remark}

\subsection{Acknowledgment}
We would like to thank Dror Bar-Natan, Vincent Franjou, Jerry Levine, 
Thang T.Q. Le, Christine Lescop, Gregor Masbaum, Paul Melvin, Xiao-Song Lin, 
Dylan
Thurston and Pierre Vogel for useful conversations.  We also wish to thank the 
referee
for numerous suggestions and comments.
\newpage
\v8
\section{Preliminaries}
\lbl{sec.pre}

\subsection{Preliminaries on Chinese characters}
\lbl{sub.cc}

Recall that a {\em Chinese character} is a  graph such that every vertex has
valency $1$ or $3$ (often called a uni-trivalent graph), together with a 
cyclic order of the edges at each of its trivalent vertices.  There is a
degree-preserving linear isomorphism 
$\chi: \B \to
\Ai$ between the graded coalgebra $\B$ of Chinese characters (modulo the
antisymmetry and IHX relations) and the graded coalgebra of chord diagrams 
$\Ai$ on a circle,
see \cite[theorem 8]{B-N}, given by mapping a Chinese character $C$ with $n$ 
legs to
$1/n!$ times the sum of the $n!$ ways of joining all of its legs to $n$ chosen
ordered points on a fixed circle.  The degree of a Chinese character or chord 
diagram is
half the number of vertices, and the primitive diagrams are the
connected ones.

Since the map $\chi$ is a vector space isomorphism, we will
{\it identify} $\B$ and  $\Ai$ via $\chi$.
Note that $\B$ has {\em two} commutative multiplications;\footnote{The
two multiplications are different.  For a conjectural  relation between these
two multiplications, see \cite[Conjecture 2]{BGRT}.}
one is induced by the multiplication on $\Ai$ via $\chi$, denoted by
$\cdot_\times$,  and the other is
the disjoint union of Chinese characters, denoted by 
$\cdot_\sqcup$. In what follows, we will suppress $\cdot_\times$ from the
notation, but will explicitly use $\cdot_\sqcup$.  Thus, 
$\exp_\sqcup$ will be used to designate the exponential with respect to the
$\cdot_\sqcup$ multiplication, and $\exp$ will be used to designate the
exponential with respect to the $\cdot_\times$ multiplication.

We will be interested in several important subspaces of $\B$.  Let $\Bi$ denote
the subspace of $\B$ which is spanned by Chinese characters, no component of 
which is
(homeomorphic to) an {\em interval}.\footnote{An interval 
is a Chinese character
of degree $1$ with $2$ univalent vertices and no trivalent ones.}  
$\Bi$ is a subalgebra of
$\B$ with respect  to
either multiplication.  Note that $\Bi$ is a direct summand of 
$\B$ with complementary factor  the span of Chinese characters which 
contain an interval component.  $\Bi$ is related to a {\em deframing} projection 
map
$F: \Ai \to \Ai$ (whose image will be denoted by $\Ap$) 
defined in \cite[part 2 of Theorem 4 and
exercise 3.16]{B-N}. Using the isomorphism $\chi$, the image of the induced 
deframing
map (also denoted by $F$) $F:\B \to \B$\footnote{Note that $F$ is {\em not} the projection in
the above  direct-sum decomposition of $\B$.} was shown in 
\cite[Corollary 4.4]{KSA} to
coincide with $\Bi$.\footnote{As an exercise, the reader may
try to find a conjectural formula for $F$ in terms of Chinese characters using 
\cite[Conjecture 2]{BGRT}.}  

Let $\Bii$ denote the subspace of $\B$ which is spanned by
Chinese characters some component of which has more trivalent vertices
than univalent ones.  $\Bii$ is a direct summand of $\Bi$.  In fact, one
has the direct sum decomposition $\Bi=\Bii\oplus \Bwh$, where
$\Bwh$\footnote{With respect to the $\cdot_\sqcup$ multiplication, $\Bwh$ is a
polynomial algebra on the  set of wheels with an even number of legs 
(the odd-legged wheels vanish by antisymmetry).} denotes the subspace of $\B$
spanned by all \cc s every component of which is a wheel (see section 2.2
below).

Let $P_{wh}: \B \to \B$ denote the composition of the deframing map $F$
followed by the projection to the subspace $\Bwh$.

We close this section with the following characterization of the algebra
of Alexander-Conway weight systems, due to \cite{KSA,Kr}.  Recall that
a {\em weight system} $W$ is a linear map $W: \Ai\to\BQ$.  Weight systems
can be multiplied and thus they form an algebra. Given a finite type
invariant of knots (or a power series of such invariants, such as the
Alexander-Conway polynomial, which will be discussed in greater detail in the
next section) there is an associated weight system, generating a subalgebra
in the algebra of weight systems. We now have the following

\begin{theorem}\cite{Kr,KSA}
\lbl{lem.con}
A weight system $W: \Ai \to \BQ$ lies in the algebra of Alexander-Conway
weight systems if and only if it factors through $P_{wh}$.
\end{theorem}

\subsection{The Alexander-Conway polynomial and its weight system}
\lbl{sub.ac}

In this section we review some well known properties of the Alexander-Conway
polynomial and its associated weight system.  (For the Alexander-Conway
polynomial and further references, see for example the exposition in the
appendix of \cite{Le}.)
The {\em Conway} polynomial $C$ \cite{Co,Ka}
of a knot (considered as a polynomial in $z$)
 is defined by the relations:
\begin{gather}
{C}\left(\printname{overcross}
	\setlength{\unitlength}{0.03\standardunitlength}
	\begin{array}{c}  \hspace{-1.7mm}
        	\raisebox{-8pt}{
\begingroup\makeatletter\ifx\SetFigFont\undefined
\def\x#1#2#3#4#5#6#7\relax{\def\x{#1#2#3#4#5#6}}%
\expandafter\x\fmtname xxxxxx\relax \def\y{splain}%
\ifx\x\y   
\gdef\SetFigFont#1#2#3{%
  \ifnum #1<17\tiny\else \ifnum #1<20\small\else
  \ifnum #1<24\normalsize\else \ifnum #1<29\large\else
  \ifnum #1<34\Large\else \ifnum #1<41\LARGE\else
     \huge\fi\fi\fi\fi\fi\fi
  \csname #3\endcsname}%
\else
\gdef\SetFigFont#1#2#3{\begingroup
  \count@#1\relax \ifnum 25<\count@\count@25\fi
  \def\x{\endgroup\@setsize\SetFigFont{#2pt}}%
  \expandafter\x
    \csname \romannumeral\the\count@ pt\expandafter\endcsname
    \csname @\romannumeral\the\count@ pt\endcsname
  \csname #3\endcsname}%
\fi
\fi\endgroup
\begin{picture}(624,639)(0,-10)
\thicklines
\path(387,237)(612,12)
\path(237,387)(12,612)
\path(118.066,548.360)(12.000,612.000)(75.640,505.934)
\path(12,12)(612,612)
\path(548.360,505.934)(612.000,612.000)(505.934,548.360)
\end{picture}
 }
        	\hspace{-1.9mm}
	\end{array}

\right)
  -{C}\left(\printname{undercross}
	\setlength{\unitlength}{0.03\standardunitlength}
	\begin{array}{c}  \hspace{-1.7mm}
        	\raisebox{-8pt}{
\begingroup\makeatletter\ifx\SetFigFont\undefined
\def\x#1#2#3#4#5#6#7\relax{\def\x{#1#2#3#4#5#6}}%
\expandafter\x\fmtname xxxxxx\relax \def\y{splain}%
\ifx\x\y   
\gdef\SetFigFont#1#2#3{%
  \ifnum #1<17\tiny\else \ifnum #1<20\small\else
  \ifnum #1<24\normalsize\else \ifnum #1<29\large\else
  \ifnum #1<34\Large\else \ifnum #1<41\LARGE\else
     \huge\fi\fi\fi\fi\fi\fi
  \csname #3\endcsname}%
\else
\gdef\SetFigFont#1#2#3{\begingroup
  \count@#1\relax \ifnum 25<\count@\count@25\fi
  \def\x{\endgroup\@setsize\SetFigFont{#2pt}}%
  \expandafter\x
    \csname \romannumeral\the\count@ pt\expandafter\endcsname
    \csname @\romannumeral\the\count@ pt\endcsname
  \csname #3\endcsname}%
\fi
\fi\endgroup
\begin{picture}(624,639)(0,-10)
\thicklines
\path(612,12)(12,612)
\path(118.066,548.360)(12.000,612.000)(75.640,505.934)
\path(387,387)(612,612)
\path(548.360,505.934)(612.000,612.000)(505.934,548.360)
\path(12,12)(237,237)
\end{picture}
 }
        	\hspace{-1.9mm}
	\end{array}

\right)
  =z {C}\left(\printname{vbridge}
	\setlength{\unitlength}{0.03\standardunitlength}
	\begin{array}{c}  \hspace{-1.7mm}
        	\raisebox{-8pt}{
\begingroup\makeatletter\ifx\SetFigFont\undefined
\def\x#1#2#3#4#5#6#7\relax{\def\x{#1#2#3#4#5#6}}%
\expandafter\x\fmtname xxxxxx\relax \def\y{splain}%
\ifx\x\y   
\gdef\SetFigFont#1#2#3{%
  \ifnum #1<17\tiny\else \ifnum #1<20\small\else
  \ifnum #1<24\normalsize\else \ifnum #1<29\large\else
  \ifnum #1<34\Large\else \ifnum #1<41\LARGE\else
     \huge\fi\fi\fi\fi\fi\fi
  \csname #3\endcsname}%
\else
\gdef\SetFigFont#1#2#3{\begingroup
  \count@#1\relax \ifnum 25<\count@\count@25\fi
  \def\x{\endgroup\@setsize\SetFigFont{#2pt}}%
  \expandafter\x
    \csname \romannumeral\the\count@ pt\expandafter\endcsname
    \csname @\romannumeral\the\count@ pt\endcsname
  \csname #3\endcsname}%
\fi
\fi\endgroup
\begin{picture}(616,631)(0,-10)
\thicklines
\path(129.845,586.698)(8.000,608.000)(105.486,531.865)
\put(-125.275,308.000){\arc{656.543}{5.1305}{7.4359}}
\path(510.514,531.865)(608.000,608.000)(486.155,586.698)
\put(741.275,308.000){\arc{656.543}{1.9889}{4.2943}}
\end{picture}
 }
        	\hspace{-1.9mm}
	\end{array}

\right),
\\ 
  C\left(\text{$c$-component unlink}\right)=\begin{cases}
    1 & \text{if }c=1 \\
    0 & \text{otherwise.}
  \end{cases} \notag
\end{gather}
With the terminology of \cite[section 3.1]{B-NG},
the Conway polynomial itself is not a
{\em canonical} Vassiliev power series, but its 
renormalized
reparametrized version 
$$
\tilde{C}(\hb)=
  \frac{\hb}{e^{\hb/2}-e^{-\hb/2}}C(e^{\hb/2}-e^{-\hb/2})
$$ is a canonical Vassiliev power series (i.e., it satisfies
$\tilde{C}(K)(h) = W_{C,\hb} \circ \ZK (K)$, see below).
Similarly, the {\em Alexander} polynomial, defined by
$A(t)=C(t^{1/2}-t^{-1/2})$, is not a canonical Vassiliev power series,
but it becomes canonical when multiplied by
$\frac{\hb}{e^{\hb/2}-e^{-\hb/2}}$ and evaluated at $t=e^\hb$
(as this product is $\tilde{C}(\hb)$).  

Let $W_C : \Ai \to \BQ$ denote the {\em weight system} of $\tilde{C}$ (which
is equal to the weight system of $C$). It has the property that it is
a deframed multiplicative weight system.  (Recall that a weight system
 $W: \Ai \to \BQ$ is called
{\em deframed} if it factors as a composition $\Ai \stackrel{F}{\to} \Ai 
\to \BQ$, where
$F$ is the  deframing map. Furthermore, a weight
system $W: \Ai \to \BQ$ is called {\em multiplicative} if
for all chord diagrams
$CD_1, CD_2$ of degrees $m_1, m_2$ respectively, we have:
$W_{m_1+ m_2}(CD_1 \cdot CD_2)=W_{m_1}(CD_1) W_{m_2}(CD_2)$.)

The weight system $W_C$ was calculated on linear chord diagrams
in \cite[Theorem 3]{B-NG}. Its 
expression in terms of Chinese characters in $\Bi$ was given by Kricker
\cite[Theorem 2.10]{Kr} as follows:
\begin{equation}
\lbl{eq.wac}
W_C (\xi) =
\begin{cases}
(-2)^p & \text{if $m=2n$ and $\xi$ is a disjoint union of} \\ 
 & \text{$p$ even-legged
wheels} \\
0 & \text{otherwise}
\end{cases}
\end{equation}
where $\omega_{2n}$ is a wheel with $2n$ legs, see Figure \ref{w2n}.

\v8
\begin{figure}[htpb]
$$ \psdraw{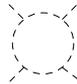}{0.4in} $$
\caption{The wheel $\omega_4$ with $4$ legs. 
Its trivalent vertices are oriented
clockwise.}\lbl{w2n}
\end{figure}
\v8

Let $\ZK: \K \to \Ai$ denote the {\em universal} \fti\ of knots, constructed by
Kontsevich
\cite{Ko} (see \cite{B-N}), and let $\xi(K) \eqdef \text{log} \ZL(K)$
denote its logarithm.  (N.b., since knots are considered as unframed, their 
image under the map $\ZK$ lies in the image of the deframing
map.\footnote{Actually, $\ZK(K)$ lies in a quotient of
$\Ai$ isomorphic to $\Ap$.}  For an extension to a 
functor $Z$ on the category of framed
$q$-tangles, see e.g.,  \cite{B2,Ca,KT,LM}. 
Then for a $0$-framed knot $K$ in $S^3$, $\ZK(K)$ coincides with the 
value  of $Z$ on $K$. In particular, $\xi(K)$ lies in (the primitive part
of, see  for example \cite {LM2}) $\Ap\simeq\Bi$.  So $P_{wh} \xi(K)$ consists
of a sum of  (even-legged) wheels.  

Define $ a_{2m}(K)$ by
\begin{equation}
\lbl{eq.axi}
\sum_{m=1}^\infty a_{2m}(K) \omega_{2m} \eqdef P_{wh} \xi(K).
\end{equation}

Let $W_{C,\hb}$ denote the product of $W_C$ and $\hb^{\deg}$,
where $\hb^{\deg}$ is the
operator that multiplies every degree $m$ diagram by $\hb^m$. 
Let $\U$ denote the (zero-framed) unknot\footnote{It may be of interest to
note that the value of $\ZK(\O)$ is  conjecturally given by the equation
$P_{wh}(\xi(\O))=\xi(\O)$.  (See \cite[Conjecture 1]{BGRT}, where it is shown 
that the conjecture holds on the level of
semisimple Lie algebras.)  Note also that the Alexander-Conway polynomial is
determined by the ${\mathfrak s \mathfrak l}_N$ colored Jones polynomial 
(see \cite{B-NG}).} 
and set $b_{2m}=a_{2m}(\U)$.  We now claim that

\v8
\begin{itemize}
\item{{\bf Fact 1.}}
For a zero framed knot in $S^3$ we have:
$$-1/2 \text{log}(A(K)(e^\hb)) = \sum_{m=1}^\infty a'_{2m}(K) \hb^{2m},$$
where $a'_{2m}(K)=a_{2m}(K)-b_{2m}$.
Indeed,  one has that
\vspace{0.3cm}\newline
\begin{tabular}{rcll} \vspace{0.1cm}
$\tilde{C}(K)(h)$ & = & $W_{C,\hb} \circ \ZK (K)$ & 
since $\tilde{C}$ is canonical  \\ \vspace{0.1cm}
& = & $\exp( W_{C,\hb} \xi(K))$ & since $W_{C,\hb}$ is multiplicative \\
& = & $\exp(-2 \sum_{m=1}^\infty a_{2m}(K)  \hb^{2m})$ & 
by equation 
\reftool.
\end{tabular}
\end{itemize}
\v8
Thus we have, 
$-1/2 \text{log}(\tilde{C}(K)(\hb)) = \sum_{m=1}^\infty a_{2m}(K) \hb^{2m}.$
   In
particular, 
$-1/2 \text{log}(\frac{\hb}{e^{\hb/2}-e^{-\hb/2}}) = \sum_{m=1}^\infty b_{2m}
\hb^{2m},$ and the result follows since 
$\tilde{C}(K)(\hb)=\tilde{C}(\U)(\hb)
A(K)(e^\hb)$.
\v8

Define $\alpha(K)$ in (the completion of) $\Bwh$ by: 
\begin{equation}
\lbl{eq.xif}
\alpha(K) =
\sum_{m=1}^\infty (b_{2m}+ a_{2m}(K)) \omega_{2m}
\end{equation}
Note that $\alpha(K)=\sum_{m=1}^\infty (2b_{2m}+ a'_{2m}(K)) \omega_{2m}.$

Similarly, let $M$ denote a 3-manifold which
satisfies $H_1(M,\BZ)=\BZ$, and let $A(M)(t)$ denote its Alexander polynomial,
normalized so that it is symmetric (in $t$ and $t^{-1}$) and evaluates to $1$ at 
$t=1$.
We define $a'_{2m}(M)$ by 
$-1/2 \text{log}(A(M)(e^\hb)) = \sum_{m=1}^\infty a'_{2m}(M)\hb^{2m}$.   

Define $\alpha(M)$ by: 
\begin{equation}
\lbl{eq.xxif}
\alpha(M)\eqdef
\sum_{m=1}^\infty (2b_{2m}+ a'_{2m}(M)) \omega_{2m}.
\end{equation}
Obviously, $\alpha(M)$ and $A(M)$ can be computed from each other.

\subsection{Preliminaries on the $LMO$ invariant}

In this section we review some well known properties of the invariant
$\ZL$.  We denote by $\{x\}_m$ the degree $m$ part of $x$. 
Recall from \cite{LMO} that for every integer $f$, and every knot
$K$ in $S^3$, the value of the universal invariant on $S^3_{K,f}$ is given
by:

$$
\ZL(S^3_{K,f})=\sum_{m=0}^{\infty} \Big\{
\frac{\i_m(c(f)\ZK(K))}{c_m(f)} \Big\}_m \in \Ao
$$ (the product in $\B$ is taken with respect to the $\cdot_\times$
multiplication), where $c(f)=\exp( \frac{f}{2}\Theta) \ZK(\O)$ and 
$c_m(f)= \i_m (\exp( \text{sgn}(f) \Theta/2 ) \ZK(\O)^2)$ (resp. 
$1$) if $ f \neq 0$ (resp. $f=0$).
Here $\Theta$ denotes the unique chord diagram of degree 1 on a
circle, $\i_m :
\Ai\to
\Ae{}$ is a map defined in 
\cite[section 2]{LMO}, $\O$
is the zero framed unknot, and $\ZK(K)$ is the value of the universal knot
invariant with the {\em zero} framing. Note that $\ZK(\O)$ is denoted by
$\nu$ in \cite{LMO,L}.

Let $\pi_m$ denote the projection $\Ao \to \Ae{\le m}$.
In the special case when $f=+1$
(n.b., the formula below holds since we are in the case of an
integral homology sphere, see \cite{LMO}), 
one has the formula
$$
\pi_m(\ZL(S^3_{K_,+1}))=\pi_m\big(\frac{\i_m(c\ZK(K))}{c_m}\big)
\in \Ae{\le m}
$$
where $c= \exp(\frac{1}{2}\Theta)\ZK(\O)$ and 
where $c_m=\i_m (\exp( \frac{1}{2}\Theta) \ZK(\O)^2)$.

\v8
The map $\i_m$, though rather complicated when
evaluated on chord diagrams on a circle, becomes more
transparent when evaluated on Chinese characters.  In particular, it follows
from its definition that for a Chinese character $C$ with $l$
legs, we have: 

\begin{equation}
\i_m (C)=
\begin{cases}
O_{-2m}(\clos{C}) & \text{if $l=2m$}\\
0 & \text{otherwise}
\end{cases}
\end{equation}
where $\clos{C}$ denotes the closure of $C$, i.e., the
sum of all $(2m-1)!!$ ways of closing its legs by joining the univalent vertices 
in pairs, and $O_{-2m}$ is the map which sets
circle components equal to $-2m$.\footnote{
To see this, see \cite{L2}, note that 
the total symmetrization of the element $T_l^m$, described 
in \cite{LMO}, vanishes, if $l$ is different from $2m$.  The formula
follows, since the total symmetrization of $T_{2m}^m$, applied to a Chinese
character, corresponds to the sum of all ways of closing up the character.}
Note that in the special case that no connected component of $C$
is an interval $I$, then 
no connected component of $\clos{C}$ is a circle, and so
$O_{-2m}(\clos{C})=\clos{C}$.  Note also that if 
$C$ has $2k$ legs, then $O_{-2m}(\clos{C\cdot_\sqcup I})=
(-2)(m-k)O_{-2m}(\clos{C})$.

We  have the following:

\begin{itemize}
\item{{\bf Fact 2.}} Fix a nonnegative integer $m$.
Given $a \in 1 + \B_{\ge 1}$, $b \in \B'_{\ge 2m}$, 
$c \in 1 + \Ae{\ge 1}$, then: 
$$\Big\{\frac{\i_m (a b))}{c}\Big\}_m  = \clos{P_{wh}(\{b\}_{2m})} 
\in \Ae{m}.$$ 
Note that this identity holds with respect to either multiplication in $\B$.
\end{itemize}
For the proof, note  that $\i_m$ reduces degree by $m$.  
In particular, the only part of $b$ which contributes to both sides 
is $\{b\}_{2m}$, the degree $2m$
part of $b$ lying in $\B'_{2m}$.
Note also that an element of $\B'_{2m}$
has at most $2m$ legs; moreover, it has exactly $2m$ legs if and only if it
lies in $(\Bwh)_{2m}$, i.e., it is a linear combination of Chinese
characters, all of whose components are even-legged wheels.  
This shows the above claim.

\section{Proofs}
\lbl{sec.pf}

In this section we give the proofs of theorems \ref{thm.al},
\ref{thm.q2} and \ref{thm.q3}.
\v8

\begin{lemma}
\lbl{lem.clos}
The map $x\mapsto \clos{x}$ from $\Bwh\to \Ae{}$ is injective, but
its restriction $(\Bwh)_{2n} \to \Ae n$ for  $n \geq 4$, is not surjective.
\end{lemma}

\begin{proof}  Note that the map $\clos{\ }$ sends connected Chinese characters
to connected trivalent graphs.  In particular it induces a map of 
primitives of 
$\Bwh$
which is easily seen to be injective.\footnote{Indeed, consider the
multiplicative  map $W :
{\Ao} \to {\BQ}[\![c]\!]$ 
defined by imposing the relation $W(H)=W(=)-W(X)$ and setting
any resulting circle
components equal to $c$.
Here $H$ denotes a diagram which 
in a neighborhood of some arc looks like an $H$, and $=$, resp. $X$, is 
obtained from $H$ by replacing this neighborhood by two arcs joining the 4 points
on the boundary of the neighborhood which are on
the same side of the arc, resp. diametrically opposed.  
Then one has that $W(\clos{w_{2}})=c^2-c$ and that $W(\clos{w_{2m+2}})=(c+2m)
W(\clos{w_{2m}})$, so $\clos{w_{2m}}\ne 0$.}

Now note further that the map $\clos{\ }$, although not multiplicative, sends a
product of primitives to the product of their closures, plus terms with fewer
connected components.  This implies that the
map $\clos{\ }$ on $\Bwh$ is injective, and further that the preimage of 
the set of primitives
is the set of primitives of $\Bwh$.  Thus, were the map also surjective as well, 
it would send the primitives onto the
primitives.  But as the primitive part of $(\Bwh)_{2n}$ is of dimension 1, 
and since the dimension of the primitive part of $\Ae{n}$ is  $>1$, for $n
\geq 4$  (see e.g. \cite{B-N}), it follows that $\clos{\ }$ is not
surjective in degree $n \geq 4$. 
\end{proof}
\v8

\begin{proof} [Proof of theorem \ref{thm.al}]

We first give the proof in case $M=S^3_{K,0}$ is obtained by surgery on
a zero-framed knot $K$ in $S^3$.  We have
\begin{eqnarray*}
\ZL(M) & = &
\sum_{m=0}^{\infty} \{ \i_m(\ZK(\O) \ZK(K)) \}_m  \\
& = & \sum_{m=0}^{\infty} \{ \i_m 
(\exp(\xi(\O))\exp(\xi(K))) \}_m  \\ 
& = & \sum_{m=0}^{\infty} \{ \i_m  (\exp(\xi(\O)+\xi(K))) \}_m  \\ 
& = & \clos{P_{wh}\exp(\xi(\O)+ \xi(K))} \\ 
& = & \clos{\exp_\sqcup(P_{wh}(\xi(\O)+\xi(K)))} \\  
& = & \clos{\exp_\sqcup(\sum_{m=1}^{\infty} (a_{2m}(\O) 
      + a_{2m}(K))w_{2m}))}\\  
& = & \clos{\exp_\sqcup\alpha(K)}
\end{eqnarray*} 
where the first and second equality is by definition, the third follows 
since we are in a commutative algebra, the fourth follows from fact 2 (with
$a=c=1$), the fifth follows since $P_{wh}$ is an algebra homomorphism, and
the last two follow from the definitions.  This shows that the
invariant $\ZL(S^3_{K,0})$ is determined by the Alexander polynomial $A(K)$.
Since $A(K)=A(M)$ (and hence $\alpha(K)=\alpha(M)$), the result follows.

Conversely, by Lemma \ref{lem.clos}, 
the map $\sum_{m=1}^{\infty} c_{2m} w_{2m}\to
\clos{\exp_\sqcup(\sum_{m=1}^{\infty} c_{2m} w_{2m})}$ is the composite of
two injective maps and hence is injective. It follows that
$\ZL(S^3_{K,0})$ determines the Alexander polynomial.

\v8
To prove the general case, first note that we may obtain
$M$ via surgery on a boundary link $K\cup L$ in $S^3$, 
where the framing on $K$ is the 
zero framing, and the framing of each component of $L$ is $\pm 1$.  
Indeed, one may obtain $M$ by zero-framed surgery on a knot in an integral
homology sphere,  which in turn may be obtained by surgery on a $\pm 1$-framed
boundary link $L$.  It suffices to  isotope, in this homology sphere, the
Seifert surface for the knot so as to be disjoint from  the Seifert surfaces of
the components of $L$.\footnote{ The same argument can be used to see that
every integral homology sphere, $\Sigma$, can be obtained by unit-framed
surgery on a boundary link:  First note that $\Sigma$  is surgery on some link
and after stabilization and handle sliding, the link may be  assumed to be $\pm
1$-framed with zero linking numbers.  In particular, $\Sigma$ can be  obtained
by a sequence of $\pm 1$-framed surgeries on knots in homology spheres. 
Arguing by  induction and applying the Seifert surface argument above,
establishes the result.}

In this case, one still has
$A(K)=A(M)$ (since the link is a boundary link, the Seifert form on $K$ in 
$S^3$ is the
same as the Seifert form of $K$ in the homology sphere obtained by surgery on 
$L$).
Moreover, since the link is boundary, its Milnor invariants vanish, and hence, 
by
\cite{HM}, $\ZK(K\cup L_0)$ consists of diagrams, none of which are trees 
(where $L_0$ denotes the link $L$ with zero framing).
Consequently, using a counting argument similar to
\cite{Ha2}, one can check that:
$\iota_m(c'(f)\ZK(K\cup L_0))=\iota_m(c'(f)\ZK(K\sqcup L_0))$
where $c'(f)$ are terms that depend on the framing of $L$, and $K \sqcup L_0$
denotes the disjoint union of the links $K$ and $L_0$. 
The definition of the $LMO$ invariant and its multiplicative
property under connected sum  implies that
$\ZL(M)=\ZL(S^3_{K,0} \sharp S^3_{L})= \ZL(S^3_{K,0})$, thus finishing 
the proof of the theorem. 
\end{proof}

\begin{proof} [Proof of theorem \ref{thm.q2}.]
Recall, \cite{L}, that $\ZL_{\le m}: \M/\M_{3m+1}\to\Ae{\le m}$ is an 
isomorphism.
We will prove that the map $\K\to\M/\M_{3m+1}$ is onto
(or equivalently, that
the composite map
$\B'\to\B'_{\le 2m}=\K/\K_{2m+1}\to\M/\M_{3m+1}=\Ae{\le m}$ is onto),
 which is dual
to the statement of theorem \ref{thm.q2}.  

The map $\B'_{\le 2m}\to \Ae{\le m}$ is given by the formula 
$$
x\mapsto \pi_m\big(\frac{\i_m(c x)}{c_m}\big)
\in \Ae{\le m},
$$
where $c,c_m$ are as in secion 2.3.

Let $x$ be a Chinese character with $2k$ legs of degree $n+k$, 
having no interval components, $n\le m$, $k\le m$. Then   
$\clos{x}$ has degree $n$.  Under the above mapping, a computation shows that
$x\mapsto (-1)^k\clos{x} + o(n+1)$, where $o(n+1)$ denotes terms of degree
$\ge n+1$.  Note that any connected graph is the closure of a
connected Chinese character with 2 legs.   Moreover, the map
$x\mapsto \clos{x}$ sends a product of connected Chinese characters (without
interval components) to the product of
their closures plus terms each of which has fewer components.
It follows by downward induction on the degree and upward induction on the 
number of
components, that the map $\B'\to \Ae{\le m}$ is surjective.
\end{proof}

\v8
\begin{proof} [Proof of theorem \ref{thm.q3}]
Consider the map $\hpsi : \K \to \Ao $ given by:
\begin{equation}
\lbl{eq.tpsi}
\hpsi(K)=  \ZL(S^3_{K,+1})
\end{equation}
as well as its truncation, $\hpsi_{\le m}=\pi_m \circ\hpsi$, where $\pi_m$
is the projection $\Ao \to \Ae{\le m}$.  Let $S$ be a $K$-admissible
$2m$-component set. 
Here $S \subset S^3$ denotes the union of $2m$ disjoint embedded balls that 
intersect the
knot in a $\pm$ crossing, and $[K,S]$ is the signed sum
of all knots obtained by changing the crossings.
(Recall that such sums generate the $2m$-th term of the Vassiliev filtration.)
One has that:
\begin{eqnarray*}
\hpsi_{\le m}([K,S])& = &
\pi_m(\sum_{S' \subseteq S} (-1)^{|S'|}\ZL(S^3_{K_{S'},+1}))  \\
& = &\pi_m \big( \frac{\i_m (c \ZK([K,S]))}{c_m} \big) \\
& = &\big\{ \frac{\i_m (c \ZK([K,S]))}{c_m} \big\}_m  \\
& = & (-1)^m \{ \clos{P_{wh}(\ZK([K,S]))} \}_m
\end{eqnarray*}
(since by \cite{LMO}, $\{c_m\}_0=(-1)^m$).

Similar computations show that $\hpsi_{\le m}([K,S])$ vanishes,
if $S$ is a $K$-admissible $n$-component set, with $n>2m$.  It follows that
$\hpsi_{\le m}$ is a $\Ae{\le m}$-valued \fti\ of knots of order $2m$.
Moreover, the above shows that the weight system $W_{\Psi,m}: \Aps{2m} \to
\Ae{m}$ factors through the projection $P_{wh}$ to $\Bwh$.
It follows from Theorem \ref{lem.con} that
$W_{\Psi,m}$ lies in the algebra of Alexander-Conway weight systems with
values in $\Ao$.
\end{proof}

\begin{remark}
\lbl{rem.gr2}
The above formula shows that the weight system is
explicitly calculated as the composite of the projection to $(\Bwh)_{2m}$
followed by the closure mapping $\clos{\ }\colon (\Bwh)_{2m}\to \Ae{m}$.
This proves the statement dual to that of Remark \ref{rem.gr}, since by Lemma
\ref{lem.clos}, 
the closure mapping on $\Bwh$ is injective, but not surjective. 
\end{remark}


\ifx\undefined\bysame
	\newcommand{\bysame}{\leavevmode\hbox 
to3em{\hrulefill}\,}
\fi

\end{document}